\documentclass{amsart}
\usepackage{head}

\begin{document}
\title{The $v_1$-Periodic Region in the cohomology of the $\C$-motivic Steenrod algebra}
\author{Ang Li}
\maketitle

\begin{abstract}
    We establish a $v_1$-periodicity theorem in $\Ext$ over the $\C$-motivic Steenrod algebra. The element $h_1$ of $\Ext$, which detects the homotopy class $\eta$ in the motivic Adams spectral sequence, is non-nilpotent and therefore generates $h_1$-towers. Our result is that, apart from these $h_1$-towers, $v_1$-periodicity operators give isomorphisms in a range near the top of the Adams chart. This result generalizes well-known classical behavior.
\end{abstract}

\section{Introduction}

\subsection{Background and Motivation}

One of the primary tools for computing stable homotopy groups of spheres is the Adams spectral sequence. The $E_2$-page of the Adams spectral sequence is given by $\Ext_{\cA^{cl}}^{*,*}(\F_2,\F_2)=H^{*,*}(\cA^{cl})$, which we denote by $\Ext_{cl}$, where $\cA^{cl}$ is the classical Steenrod algebra. For $\Ext_{cl}$, Adams \cite{Per} showed that there is a vanishing line of slope $\frac{1}{2}$ and intercept $\frac{3}{2}$, and J. P. May showed there is a periodicity line of slope $\frac{1}{5}$ and intercept $\frac{12}{5}$, where the periodicity operation is defined by the Massey product $P_r(-):=\langle h_{r+1}, h_0^{2^r},-\rangle$. This result has not been published by May, but can be found in the thesis of Krause:

\begin{thm}\label{cpt}\cite[Theorem 5.14]{Kra}
For $r\geq 2$, the Massey product operation $P_r(-):=\langle h_{r+1}, h_0^{2^r},-\rangle$ is uniquely defined on $\Ext_{cl}^{s,f}=H^{s,f}(\cA^{cl})$ when $s>0$ and $f>\frac{1}{2}s+3-2^r$, where $s$ is the stem, and $f$ is the Adams filtration.

Furthermore, for $f>\frac{1}{5}s+\frac{12}{5}$, the operation 
\[P_r\colon H^{s,f}(\cA^{cl})\xrightarrow{\iso} H^{s+2^{r+1},f+2^r}(\cA^{cl})\]
is an isomorphism.
\end{thm}

The purpose of this article is to discuss an analog of the theorem above in the $\C$-motivic context. Motivic homotopy theory, also known as $\A^1$-homotopy theory, is a way to apply the techniques of algebraic topology, specifically homotopy, to algebraic varieties and, more generally, to schemes. The theory was formulated by Morel and Voevodsky \cite{MV}.\\

In this paper we analyze the case where the base field $F$ is the complex numbers $\C$. Let $\M_2$ denote the bigraded motivic cohomology ring of Spec $\C$, with $\F_2=\Z/2$-coefficients. Voevodsky \cite{Voe} proved that $\M_2\iso \F_2[\tau]$. Let $\cA$ be the mod 2 motivic Steenrod algebra over $\C$. The motivic Adams spectral sequence is a trigraded spectral sequence with \[E_2^{*,*,*}=\Ext_\cA^{*,*,*}(\M_2,\M_2),\] where the third grading is the motivic weight. (See Dugger and Isaksen \cite{DI1}). The $\C$-motivic $E_2$-page, which we denote by $\Ext$, has a vanishing line computed by Guillou and Isaksen \cite{GI1}. Quigley has a partial result for the motivic periodicity theorem in the case $r=2$ \cite[Corollary 5.4]{JD}.\\

The multiplication by 2 map $S^{0,0}\xrightarrow{2}S^{0,0}$ is detected by $h_0$, and the Hopf map $S^{1,1}\xrightarrow{\eta}S^{0,0}$ is detected by $h_1$ in $\Ext$. These elements have degrees $(0,1,0)$ and $(1,1,1)$ respectively. By an infinite $h_1$-tower we will mean a non-zero sequence of elements of the form $h_1^kx$ in $\Ext$ with $k\geq 0$, where $x$ is not $h_1$ divisible. We will write $h_1$-towers for infinite $h_1$-towers, and refer to $x$ as the base of the $h_1$-tower $h_1^kx$ ($k\geq 0$). A short discussion on the $h_1$-towers can be found in subsection \ref{future}. Since all $h_1$-towers are $\tau$-torsion, one might guess that the motivic $\Ext$ groups differ from the classical $\Ext^{cl}$ groups by only infinite $h_1$-towers. This is not true, but we may expect the $h_1$-torsion part of $\Ext$ to obtain a pattern similar to $\Ext^{cl}$. Our result pertains solely to this $h_1$-torsion region. 

\begin{rmk}
Let $\cA\dual$ denote the dual Steenrod algebra. For $\Ext$, we can work over $\cA\dual$ instead of $\cA$. i.e. \[E_2^{*,*,*}\iso \Ext_{\cA\dual}^{*,*,*}(\M_2,\M_2)\dual.\]
Here we view $\M_2$ as the homology of the motivic sphere instead of the cohomology; this is an $\cA\dual$-comodule.\\
\end{rmk}

The goal of this paper is the following theorem:
\begin{thm}\label{mpt}
For $r\geq 2$, the Massey product operation $P_r(-):=\langle h_{r+1}, h_0^{2^r},-\rangle$ is uniquely defined on $\Ext^{s,f,w}=H^{s,f,w}(\cA)$ when $s>0$ and $f>\frac{1}{2}s+3-2^r$.

Furthermore, for $f>\frac{1}{5}s+\frac{12}{5}$, the restriction of $P_r$ to the $h_1$-torsion 
\[P_r\colon [H^{s,f,w}(\cA)]_{h_1-\text{torsion}}\to [H^{s+2^{r+1},f+2^r,w+2^r}(\cA)]_{h_1-\text{torsion}}\]
is an isomorphism.
\end{thm}

There are close connections between the classical Adams spectral sequence and the motivic Adams spectral sequence. For instance, by inverting $\tau$ in $\Ext$, we obtain $\Ext^{cl}$. There are also abundant connections between the $\C$-motivic $\Ext$ groups, the $\R$-motivic $\Ext$ groups and the $C_2$-equivariant $\Ext$ groups. The $\rho$-Bockstein spectral sequence \cite{MH} takes the $\C$-motivic $\Ext$ groups as input and computes the $\R$-motivic $\Ext$ groups. The $C_2$-equivariant $\Ext$ groups can then be obtained \cite{GHIR} by calculating $\R$-motivic $\Ext$ groups for a negative cone. Our periodicity results ought to be relevant for future computations in $\R$-motivic and $C_2$-equivariant homotopy theory.\\

\subsection{Further Considerations}\label{future}

We study the $h_1$-torsion part of $Ext$; the $h_1$-periodic part has been entirely computed in \cite{GI2}. 
\begin{thm}\cite[Theorem 1.1]{GI2}\label{h1inv}
The $h_1$-inverted algebra $\Ext_\cA[h_1\inv]$ is a polynomial algebra over $\F_2[h_1^{\pm1}]$ on generators $v_1^4$ and $v_n$ for $n\geq 2$, where:
\begin{enumerate}
    \item $v_1^4$ is in the $8$-stem and has Adams filtration $4$ and weight $4$.
    \item $v_n$ is in the $(2^{n+1}-2)$-stem and has Adams filtration $1$ and weight $2^n-1$.
\end{enumerate}
\end{thm}

It is straightforward that $P_r$ acts injectively on the $h_1$-inverted $\Ext$; that is, $P_r$ sends an $h_1$-tower $h_1^kx$ ($k\geq 0$) to another $h_1$-tower $h_1^ly$ ($l\geq 0$). But the base $x$ might not be sent to the base $y$. As for surjectivity, there are $h_1$-towers not in the image of $P_r$, such as the $h_1$-tower on $c_0$; those are not multiples of $v_1^4$ in the $h_1$-inverted $\Ext$. Partial results about the bases of those $h_1$-towers can be found in \cite{HT}. We expect that the determination of the bases of the $h_1$-towers will lead to a complete understanding of the region in which the $v_1$-periodicity operator acts as an isomorphism on $Ext$.\\

There is another periodicity element $w_1$ in motivic $Ext$, which does not exist classically. Analogously to the Massey product $P_2(-):=\langle h_3, h_0^4,-\rangle$, there is another Massey product $g(-):=\langle h_4,h_1^4,-\rangle$. For many values of $x$, $P_2(x)$ is detected by $Px$, where $P=h_{20}^4$ has degree $(8,4,4)$ in the May spectral sequence. Similarly, for many values of $x$, $g(x)$ will be detected by $h_{21}^4\cdot x$, where $h_{21}^4$ has degree $(20,4,12)$ in the May spectral sequence. The obstruction to studying $w_1$-periodicity is that $g$ has a relatively low slope. Thus the method in this paper is not applicable. In addition, our method relies on a computation involving $\Ext_{\cA(1)\dual}$, but $g$ restricts to zero in that group. Thus a strategy for studying $g$-periodicity would need to begin with $\Ext_{\cA(2)\dual}$, which is
much more complicated \cite{Isa}.

\subsection{Organization}

We follow the approach of \cite{Kra} primarily. In Section 2, we briefly introduce the stable (co)module category, in which we can consider the $h_0$ or $h_1$-torsion part of $\Ext$ by taking sequential colimits. In Section 3, we establish the existence of a homological self-map $\theta$ and use this to show that $P_r(-)$ is uniquely defined. In Section 4, we explicitly show where $\theta$ is an isomorphism over $\cA(1)\dual$, and obtain a region where it is an isomorphism over $\cA\dual$ by moving along the Cartan-Eilenberg spectral sequence. In Section 5, we combine the results of the previous two sections together to get the motivic periodicity theorem \ref{mpt}.

\subsection{Acknowledgements}

The author would like to thank Bertrand Guillou for useful instructions and helpful discussions. The author also benefited from discussions with J.D. Quigley, Eva Belmont, and Prasit Bhattacharya and appreciates their assistance. The author thanks the referee for providing detailed comments that helped to improve the exposition of the article.

\section{Working environment: the stable (co)module category $Stab(\Gamma)$}

In order to restrict to working with only the $h_1$-torsion (also $h_0$-torsion) part, first we would like to choose a suitable working environment: a category with some nice properties that will serve our purposes. Usually $\Ext$ is defined in the derived category of $\cA\dual$-comodules, which we denote $D(\cA\dual)$. However, the coefficient ring $\M_2$ is not compact in $D(\cA\dual)$, which means that $\M_2$ does not interact well with colimits. The stable comodule category will better serve our purposes. That is a category $\sC$ such that:
\begin{enumerate}
    \item If $M$ is a $\cA\dual$-comodule that is free of finite rank over $\M_2$ and $N$ is a $\cA\dual$-comodule, then $\Hom_\sC(M,N)\iso\Ext_{\cA\dual}(M,N)$. 
    \item If $M$ is a $\cA\dual$-comodule that is free of finite rank over $\M_2$, then $M$ is compact in $\sC$. That is to say, for any sequential colimit in $\sC$ of $\cA\dual$-comodules
    \[\underset{i}{\colim} N_i:=\colim (N_0\xrightarrow{f_0}N_1\to \cdots\to N_i\xrightarrow{f_i}\cdots),\] we have $\underset{i}{\colim}\Ext_{\cA\dual}(M,N_i)\iso\Hom_\sC(M,\underset{i}{\colim} N_i)$
\end{enumerate}

The correct choice of $\sC$ is called $Stab(\cA\dual)$. The category can be constructed in various ways (see \cite[Sec. 2.1]{Bel} for details), and has several useful properties for our case. The following proposition summarizes some of the discussion in \cite[Sec. 4]{BHV}:

\begin{prop}
\label{smc}
The category $Stab(\cA\dual)$ satisfies conditions $(1)$ and $(2)$ above. 
\end{prop}

Namely, for a Hopf algebra $\Gamma$ and comodule $M$ that is free of finite rank, we have a diagram 
\[\xymatrix{
D(\Gamma)\ar@/_1pc/[dr]_-{\Hom_{D(\Gamma)}(iM,-)}& Comod_\Gamma\ar[l]_i\ar@{.>}[r]^j\ar[d]_{\Ext_\Gamma(M,-)}& Stab(\Gamma)\ar@/^1pc/[dl]^-{\Hom_{Stab(\Gamma)}(jM,-)}\\
  &\mathbf{grAb}&  \\
}
\]
where $i$ is the canonical functor and $j$ is well-defined only for comodules that are free of finite rank over $\M_2$. This diagram commutes. Because the stable comodule category cooperates nicely with taking colimits, we can compute the colimit of a sequence of $\Ext_\Gamma(M,N)$.\\

Here we introduce notation that will be used in future sections.

\begin{notn}\label{notation}
For a spectrum $M$ such that $H_*(M)$ is free of finite rank over $\M_2$, let $M$ also denote the embedded image of the homology of the spectrum $M$ in the stable comodule category (i.e., $M=j(H_*(M))$). We use $[M,N]_{*,*,*}^\Gamma$ to denote $\Hom_{Stab(\Gamma)}(M,N)$, where $M$, $N\in Stab(\Gamma)$. For example, if $M=S^0$, then $H_*(S^0)=\M_2$, which we also denote by $S$. Thus $\Ext_{\cA\dual}^{s,f,w}(\M_2,\M_2) = [S,S]^{\cA\dual}_{s,f,w}$.
When $\Gamma$ is the motivic dual Steenrod algebra, we omit the superscript $\Gamma$. This notation is consistent with \cite{Kra}.\\

We use the grading $(s,f,w)$, where $s$ is the stem, $f$ is the Adams filtration and $w$ is the motivic weight. Notice that $t=s+f$ is the internal degree. Given a self-map $\theta$: $\Sigma^{s_0,f_0,w_0}M\xrightarrow{\theta} M$ in $Stab(\cA\dual)$, we have a cofiber sequence $\Sigma^{s_0,f_0,w_0}M\xrightarrow{\theta} M\to M/\theta$ in $Stab(\cA\dual)$. The associated long exact sequence will be indexed as follows:
\[
\cdots\to [M,N]_{s+s_0+1,f+f_0-1,w+w_0}\to [M/\theta,N]_{s,f,w}\to [M,N]_{s,f,w}\to [M,N]_{s+s_0,f+f_0,w+w_0}\to \cdots
\]
Sometimes we omit indices when there is no risk of confusion. 
\end{notn}

\section{Self-maps and Massey products}
In this section, we show that the cofiber $S/h_0^k$ admits a self-map and identify it with the Massey product in Theorem~\ref{mpt}. Self-maps are maps of suspensions of an object to itself. For a dualizable object $Y$, self maps $\Sigma^nY\to Y$ can also be described as elements of $\pi_*(Y\otimes DY)$, with $DY$ the $\otimes$-dual of $Y$. In this paper we mainly deal with homological self-maps in $Stab(\cA\dual)$.\\

When considering the vanishing region and the periodicity region, we only work with the $h_0$-torsion part. (Of course, this is not much of a loss: as classically, the only $h_0$-local elements are in the 0-stem.) We next investigate the $h_1$-torsion part inside the $h_0$-torsion. For this purpose, we introduce the following notion.

\begin{defn}\label{F0}
Let $F_0$ be the fiber of $S\to S[h_0\inv]$, where $S[h_0\inv]:=\colim(S^0\xrightarrow{h_0}S\xrightarrow{h_0}\cdots)$ in $Stab(\cA\dual)$. Similarly, let $F_{01}$ be the fiber of $F_0\to F_0[h_1\inv]$ with $F_0[h_1\inv]$ defined as an analogous colimit.
\end{defn}

The group $[S,F_{01}]$ contains the subset of $[S,S]$ consisting of elements that are both $h_0$- and $h_1$-torsion, as well as the negative parts of those $h_0$ and $h_1$-towers in $F_0[h_1\inv]$. The regions we are considering are unaffected. We display the corresponding $\Ext$ groups in Figure \ref{sf0} and \ref{sf01}.
\\

\begin{figure}[h]
\centering
\hspace{-1cm}
\begin{minipage}{0.4\linewidth}
\begin{tikzpicture}[scale=0.4]

\draw[thick,->](-2,0)--(13,0);
\draw[thick,->](0,-3)--(0,8);
\draw[<->] (13,0) node[below]{$s$} (5,0) node[below]{5} (10,0) node[below]{10} (0,8) node[left]{$f$}; 
\foreach \t in{-1,0,...,12}\draw(\t,0.1)--(\t,-0.1);
\foreach \t in{-2,-1,0,...,7}\draw(0.1,\t)--(-0.1,\t);

\draw[fill=black,black](1,1)circle(0.15);
\draw[fill=black,black](2,2)circle(0.15);
\draw[fill=black,black](3,3)circle(0.15);
\draw[fill=black,black](3,1)circle(0.15);
\draw[fill=black,black](3,2)circle(0.15);
\draw[fill=black,black](6,2)circle(0.15);
\draw[fill=black,black](7,1)circle(0.15);
\draw[fill=black,black](7,2)circle(0.15);
\draw[fill=black,black](7,3)circle(0.15);
\draw[fill=black,black](7,4)circle(0.15);
\draw[fill=black,black](8,2)circle(0.15);
\draw[fill=black,black](8,3)circle(0.15);
\draw[fill=black,black](9,3)circle(0.15);
\draw[fill=black,black](9,4)circle(0.15);
\draw[fill=black,black](9,5)circle(0.15);
\draw[fill=black,black](10,6)circle(0.15);
\draw[fill=black,black](11,5)circle(0.15);
\draw[fill=black,black](11,6)circle(0.15);
\draw[fill=black,black](11,7)circle(0.15);

\draw[ultra thick,red,->](3,3)--(4,4);
\draw[ultra thick,red,->](9,4)--(10,5);
\draw[ultra thick,red,->](11,7)--(12,8);

\draw[ultra thick,black](1,1)--(3,3);
\draw[ultra thick,black](3,1)--(3,3);
\draw[ultra thick,black](3,1)--(9,3);
\draw[ultra thick,black](7,1)--(7,4);
\draw[ultra thick,black](7,1)--(9,3);
\draw[ultra thick,black](8,3)--(9,4);
\draw[ultra thick,black](9,5)--(11,7);
\draw[ultra thick,black](11,5)--(11,7);

\draw[ultra thick,blue,->](-1,0)--(-1,-3);

\end{tikzpicture}

\caption{$[S,F_0]^{\cA^\vee}_{*,*,*}$}
\label{sf0}
\end{minipage}
\hspace{1cm}
\begin{minipage}{0.4\linewidth}
\begin{tikzpicture}[scale=0.4]

\draw[thick,->](-2,0)--(13,0);
\draw[thick,->](0,-3)--(0,8);
\draw[<->] (13,0) node[below]{$s$} (5,0) node[below]{5} (10,0) node[below]{10} (0,8) node[left]{$f$}; 
\foreach \t in{-1,0,...,12}\draw(\t,0.1)--(\t,-0.1);
\foreach \t in{-2,-1,0,...,7}\draw(0.1,\t)--(-0.1,\t);

\draw[fill=black,black](1,1)circle(0.15);
\draw[fill=black,black](2,2)circle(0.15);
\draw[fill=black,black](9,4)circle(0.15);
\draw[fill=black,black](10,6)circle(0.15);
\draw[fill=black,black](3,3)circle(0.15);
\draw[fill=black,black](3,1)circle(0.15);
\draw[fill=black,black](3,2)circle(0.15);
\draw[fill=black,black](6,2)circle(0.15);
\draw[fill=black,black](7,1)circle(0.15);
\draw[fill=black,black](7,2)circle(0.15);
\draw[fill=black,black](7,3)circle(0.15);
\draw[fill=black,black](7,4)circle(0.15);
\draw[fill=black,black](8,2)circle(0.15);
\draw[fill=black,black](8,3)circle(0.15);
\draw[fill=black,black](9,3)circle(0.15);
\draw[fill=black,black](9,5)circle(0.15);
\draw[fill=black,black](11,5)circle(0.15);
\draw[fill=black,black](11,6)circle(0.15);
\draw[fill=black,black](11,7)circle(0.15);

\draw[ultra thick,black](3,1)--(3,2);
\draw[ultra thick,black](3,1)--(9,3);
\draw[ultra thick,black](7,1)--(7,4);
\draw[ultra thick,black](7,1)--(9,3);
\draw[ultra thick,black](11,5)--(11,6);
\draw[ultra thick,black](1,1)--(3,3);
\draw[ultra thick,black](8,3)--(9,4);
\draw[ultra thick,black](9,5)--(11,7);

\draw[ultra thick,black](3,2)--(3,3);
\draw[ultra thick,black](11,6)--(11,7);

\draw[ultra thick,blue,->](-1,0)--(-1,-3);
\draw[ultra thick,pink,->](-1,1)--(-2,0);
\draw[ultra thick,pink,->](6,3)--(5,2);
\draw[ultra thick,pink,->](7,5)--(6,4);

\end{tikzpicture}

\caption{$[S,F_{01}]^{\cA^\vee}_{*,*,*}$}
\label{sf01}
\end{minipage}

\end{figure}

The periodicity operator $P$ corresponds to multiplying by the element $h_{20}^4$ of the May spectral sequence, meaning that for many values of $x$, $h_{20}^4x\in \langle h_3,h_0^4,x\rangle$. However, $h_{20}^4$ does not survive to $\Ext$. As a result, multiplying by $P$ is not a map from $[S,S]$ to $[S,S]$. Luckily, \cite[Figure 2]{GI1} shows that $P$ survives in $[S/h_0,S]$. Similarly, we have the following proposition:
\begin{prop}\label{exist}
The element $h_{20}^{2^r}$ survives the May spectral sequence to $[S/h_0^k,S]$ for $k\leq 2^r$, and thus gives a corresponding element $P^{2^{r-2}}$ in $[S/h_0^k,S/h_0^k]$, i.e. a self-map of $S/h_0^k$.
\end{prop}

If $N$ is an $\cA\dual$-comodule in $Stab(\cA\dual)$, then $[S/h_0^k,S/h_0^k]$ acts on $[S/h_0^k,N]$. The corresponding element $P$ (or some power of $P$) inside $[S/h_0^k,S/h_0^k]$ induces a map from $[S/h_0^k,N]$ to itself. We would like to show that for any $k\leq 2^r$ and  $r\geq 2$, multiplying by $P^{2^{r-2}}$ on $[S/h_0^k, S]$ coincides with the Massey product $P_r(-):=\langle h_{r+1}, h_0^{2^r},-\rangle$ in a certain region. In other words, we must show that there is zero indeterminacy.

The Massey product is defined on the kernel
of $h_0^{2^r}$ on $[S,S]$, which we will denote $ker(h_0^{2^r})$. It lands in the cokernel of multiplication by $h_{r+1}$:
\[P_r(-): ker(h_0^{2^r})\to [S, S]/h_{r+1.}\]

\begin{rmk}
Originally one would like to consider the following square and see that it commutes in a certain region
\[\xymatrix{[S/h_0^k,S]\ar[r]^{-\cdot P^{2^{r-2}}}\ar[d]&[S/h_0^k,S]\ar[d]\\
ker(h_0^{2^r})\ar[r]_-{P_r(-)}&[S, S]/h_{r+1}.}
\]
The vertical maps are induced by $S\to S/h_0^k$. However, since we lost the advantage of a vanishing region of $f>\frac{1}{2}s+\frac{3}{2}$ that we need in the classical setting, the region where the vertical maps are isomorphisms is not satisfactory. We solve this problem by restricting attention to the $h_0$ and $h_1$-torsion.
\end{rmk}

To better fit our purposes, consider the Massey product defined on $[S,F_{01}]$\[P_r(-): ker_{F_{01}}(h_0^{2^r})\to [S, F_{01}]/h_{r+1.}\]
This gives the following squares, over which we have more control:

\begin{equation}
\xymatrix{[S/h_0^k,F_{01}]\ar[r]^{-\cdot P^{2^{r-2}}}\ar[d]&[S/h_0^k,F_{01}]\ar[d]\\
ker_{F_{01}}(h_0^{2^r})\ar[r]_-{P_r(-)}\ar[d]&[S, F_{01}]/h_{r+1}\ar[d]\\
ker_S(h_0^{2^r})\ar[r]_-{P_r(-)}&[S, S]/h_{r+1}}
\label{topsquare}
\end{equation}

The canonical map $F_{01}\to S$ induces a map $[S,F_{01}]\to[S,S]$ given by inclusion on the $h_0-$ and $h_1$-torsion elements and which sends negative towers to zero. The bottom square commutes for $s>0$ and $f>0$ modulo potential indeterminacy. We would like to show that the indeterminacy vanishes under some conditions.

Let $C(\eta)$ denote the cofiber of the first Hopf map \[S^{1,1}\xrightarrow{\eta}S^{0,0}.\] Writing $C_\eta$ for the cohomology $H^{*,*}(C(\eta))$, we have the following result:

\begin{thm}\cite[Theorem 1.1]{GI1}\label{vceta}
The group $\Ext_\cA^{s,f,w}(\M_2,C_\eta)$ vanishes when $s>0$ and $f>\frac{1}{2}s+\frac{3}{2}$.
\end{thm}

Theorem \ref{vceta} gives us that $[S, C_\eta]_{s,f,w}$ vanishes when $s>0$ and $f>\frac{1}{2}s+\frac{3}{2}$. In other words, there are only $h_1$-towers when $s>0$ and $f>\frac{1}{2}s+\frac{3}{2}$ in $[S,S]_{s,f,w}$. Moreover, we have the following fact:

\begin{prop}[Corollary of {\cite[Theorem 1.1]{GI2}}]
For $r\geq 1$, $h_{r+1}$ does not support an $h_1$-tower.
\end{prop}

Therefore the indeterminacy $(h_{r+1}[S,S])_{s,f,w}$ must vanish when $f>\frac{1}{2}s+3-2^r$, under the following two conditions: that $h_{r+1}$ has $s=2^{r+1}-1$, and that there are only $h_1$-towers in $[S,S]_{s,f,w}$ when $s>0$ and $f>\frac{1}{2}s+\frac{3}{2}$, which are $h_{r+1}$-torsion groups.

\begin{rmk}\label{ind}
It is easy to see that the indeterminacy $(h_{r+1}[S,F_{01}])_{s,f,w}$ also vanishes when $f>\frac{1}{2}s+3-2^r$.
\end{rmk}

The first row of the top square in \eqref{topsquare} is multiplication by some power of the element $P$. We next determine when the vertical maps are isomorphisms. 

\begin{lemma}[Motivic version of {\cite[Lemma 5.2]{Kra}}]\label{lemma 5.2}
Let $M,N\in Stab(\cA\dual)$.
Assume that $[M,N]$ vanishes when $f>as+bw+c$ for some $a,b,c\in \R$, let $\theta:\Sigma^{s_0,f_0,w_0}M\to M$ be a map with $f_0>as_0+bw_0$, and let $M/\theta$ denote the cofiber of $\Sigma^{s_0,f_0,w_0}M\xrightarrow{\theta} M$.
Then \[[M/\theta,N]\to [M,N]\]
is an isomorphism above a vanishing plane parallel with the one in $[M,N]$ but with $f$-intercept given by $c-(f_0-as_0-bw_0)$.
\end{lemma}

\begin{pf}
The result follows from the long exact sequence associated to the cofiber sequence $\Sigma^{s_0,f_0,w_0}M\xrightarrow{\theta} M\to M/\theta$:
\[
\cdots\to [M,N]_{s+s_0+1,f+f_0-1,w+w_0}\to [M/\theta,N]_{s,f,w}\to [M,N]_{s,f,w}\to [M,N]_{s+s_0,f+f_0,w+w_0}\to \cdots
\]
\end{pf}

\begin{rmk}
This approach could also apply to a vanishing region above several planes or even a surface. The vanishing condition of Lemma \ref{lemma 5.2} could be rephrased as the following:

Assume that $[M,N]_{*,*,*}$ vanishes when $f>\varphi(s,w)$ where $\varphi:\R^2\to\R$ is a smooth function. Then the gradient $v(-,-)=(\frac{\bdry\varphi}{\bdry s}(-),\frac{\bdry\varphi}{\bdry w}(-))$ is a vector field. Let $d=\underset{(s_0,w_0)}{max}|v(s_0,w_0)|$, and assume both $\frac{f_0}{s_0}$ and $\frac{f_0}{w_0}$ are greater than $d$. The remaining proof would follow similarly, with the $f$-intercept given by $max \{c-(f_0-ds_0),c-(f_0-dw_0)\}$.
\end{rmk}

We have this as a corollary:
\begin{cor}[Motivic version of {\cite[Lemma 5.9]{Kra}}]\label{ud}
Let $k\geq 1$. For $f>\frac{1}{2}s+\frac{3}{2}-k$, the natural map $[S/h_0^k,F_{01}]_{s,f,w}\to [S,F_{01}]_{s,f,w}$ is an isomorphism.
\end{cor}

\begin{pf}
To determine this, we need to confirm that $[S,F_{01}]_{s,f,w}$ admits a vanishing region of $f>\frac{1}{2}s+\frac{3}{2}$. The fiber sequence $F_{01}\to F_0\hookrightarrow F_0[h_1\inv]$ gives us an exact sequence: \[\cdots\to [S,F_{01}]_{s,f,w}\to [S, F_0]_{s,f,w}\xhookrightarrow{h_1\inv}[S,F_0[h_1\inv]]_{s,f,w}\to[S, \Sigma^{1,-1,0}F_{01}]_{s,f,w}\to \cdots\]
Since $[S,F_0]$ differs from $[S,S]$ only in the first column, there are only $h_1$-towers when $f>\frac{1}{2}s+\frac{3}{2}$. And by Theorem \ref{vceta} again, $[S,C_\eta]_{s,f,w}$ vanishes when $s>0$ and $f>\frac{1}{2}s+\frac{3}{2}$. In other words, above the plane $f=\frac{1}{2}s+\frac{3}{2}$, multiplying by $h_1$, which detects $\eta$, is an isomorphism from $[S, F_0]_{s,f,w}$ to $[S,F_0]_{s+1,f+1,w+1}$.

As a result, inverting $h_1$ would be an isomorphism from $[S,F_0]_{s,f,w}$ to $[S,F_0[h_1\inv]]_{s,f,w}$ when $f>\frac{1}{2}s+\frac{3}{2}$. Therefore, $[S,F_{01}]_{s,f,w}$ vanishes when $f>\frac{1}{2}s+\frac{3}{2}$. Applying Lemma \ref{lemma 5.2} gives the corollary.
\end{pf}

The results in \ref{exist} and \ref{ind} locate the region where both squares commute, thus obtaining the first part of Theorem \ref{mpt}.

\begin{thm}[Motivic version of {\cite[Proposition 5.12]{Kra}}]\label{selfmap}
For $k\leq 2^r$ and $r\geq 2$, the cofiber $S/h_0^k$ admits a self-map $P^{2^{r-2}}$ of degree $(2^{r+1}, 2^r, 2^r)$. Thus, for any $N\in Stab(\cA\dual)$, composition with $P^{2^{r-2}}$ defines a self-map on $[S/h_0^k,N]$.

When $f>\frac{1}{2}s+3-k$, the induced map coincides with the Massey product $P_r(-):=\langle h_{r+1}, h_0^{2^r},-\rangle$ with zero indeterminacy.
\end{thm}

\section{Colimits and the Cartan-Eilenberg spectral sequence}

We will obtain a vanishing region for $[S/(h_0,P),F_{01}]_{*,*,*}$ in this section. Consider the colimit \[F_0/h_1^\infty:=\underset{i}{\colim}(\Sigma^{-1,-1,-1}F_0/h_1\xrightarrow{h_1}\cdots\xrightarrow{h_1}\Sigma^{-i,-i,-i}F_0/h_1^i\xrightarrow{h_1}\cdots)\] in $Stab(\cA\dual)$. As we show in the following result, it differs from $F_{01}$ by a suspension in the region we are considering.
\begin{prop}\label{shift}
When $f>\frac{1}{2}s+\frac{3}{2}$,
\[[S,\Sigma^{-1,1,0}F_0/h_1^\infty]_{s,f,w}\iso [S,F_{01}]_{s,f,w}\]
\end{prop}
\begin{pf}
To see this, note that the colimit $F_0/h_1^\infty$ is a union of all the $h_1$-torsion in $F_0$, while the fiber $F_{01}$ detects the $h_1$-torsion together with those negative $h_1$-towers.
\end{pf}

Note that $F_0$ coincides with \[\Sigma^{-1,1,0}S/h_0^\infty:=\Sigma^{-1,1,0} \underset{i}{\colim}(\Sigma^{0,1,0}S/h_0\xrightarrow{h_0}\cdots\xrightarrow{h_0}\Sigma^{0,i,0}S/h_0^i\xrightarrow{h_0}\cdots),\] if we ignore the negative $h_0$-tower. That is, we have $[S,\Sigma^{-1,1,0} S/h_0^\infty]_{s,f,w}\iso [S,F_0]_{s,f,w}$ when $f>0$.

\begin{rmk}
We have shown that the map $[S/h_0^k,F_0/h_1^\infty]_{s,f,w}\to[S,F_0/h_1^\infty]_{s,f,w}$ is an isomorphism when $f>\frac{1}{2}s+3-k$. We consider this colimit because it is better for computational purposes (the fiber $F_{01}$ is harder to deal with than the colimit $F_0/h_1^\infty$).
\end{rmk}

Let $\theta$ be a self-map of $S/h_0^k$, and consider the cofiber sequence $S/h_0^k\xrightarrow{\theta}S/h_0^k\to S/(h_0^k,\theta)$. The vanishing region for $[S/(h_0^k,\theta),F_0/h_1^\infty]_{*,*,*}$ is the region where \[[S/h_0^k,F_0/h_1^\infty]_{s,f,w}\xrightarrow{\theta}[S/h_0^k,F_0/h_1^\infty]_{s+s_0,f+f_0,w+w_0}\] is an isomorphism. The goal of this section is to obtain a vanishing region for $[S/(h_0^k,\theta),F_0/h_1^\infty]_{*,*,*}$ in the case $k=1$ and $\theta=P$.\\

The dual Steenrod algebra is too large to work with, so we would like to start with a smaller one, namely $\cA(1)\dual\iso \M_2[\tau_0,\tau_1,\xi_1]/(\tau_0^2=\tau\xi_1,\tau_1^2,\xi_1^2)$. Then for $\cA\dual$-comodules $M$ and $N$ (thus also $\cA(1)\dual$-comodules), we can recover $[M,N]^{\cA\dual}$ from $[M,N]^{\cA(1)\dual}$ via infinitely many Cartan-Eilenberg spectral sequences along normal extensions of Hopf algebras.

A brief introduction to the Cartan-Eilenberg spectral sequence (see \cite[Ch.XV]{CE} for details) is relevant at this point. Given an extension of Hopf algebras over $\M_2$
\[E\to \Gamma \to C\]
(so in particular $E\iso \Gamma \square_C \M_2$), 
the Cartan-Eilenberg spectral sequence computes $Cotor_\Gamma(M,N)$ for a $\Gamma$-comodule $M$ and an $E$-comodule $N$. The spectral sequence arises from the double complex ($\Gamma$-resolution of $M$)$\square_\Gamma$($E$-resolution of $N$), and we have $Cotor_\Gamma(M,N)\iso \Ext_\Gamma(M,N)$ when $M$ and $N$ are $\tau$-free. 

The Cartan-Eilenberg spectral sequence has the form 
\[E_1^{s,t,*,*}=Cotor_C^{t,*}(M,\Bar{E}^{\otimes s}\otimes N)\Rightarrow Cotor_\Gamma^{s+t,*}(M,N).\]
If $E$ has trivial $C$-coaction, then we have $E_1^{s,t,*,*}\iso Cotor_C^{t,*}(M,N)\otimes \Bar{E}^{\otimes s}$. Taking the cohomology we obtain the $E_2$-page: \[E_2^{s,t,*,*}=Cotor_E^{s,*}(\M_2,Cotor_C^{t,*}(M,N))\iso \Ext_E^{s,*}(\M_2,\M_2)\otimes \Ext_C^{t,*}(M,N).\]

Let $N=F_0/h_1^\infty$. We will compute $[S/h_0,F_0/h_1^\infty]^{\cA(1)\dual}$ as an intermediate step before reaching our goal of $[S/(h_0,P),F_0/h_1^\infty]^{\cA(1)\dual}$. As a starting point, we can compute $[S/h_0,F_0]$ over $\cA(1)\dual$, via the cofiber sequence $S\xrightarrow{h_0}S\to S/h_0$.

\begin{figure}[h]
\centering
\begin{tikzpicture}[scale=0.6]

\draw[thick,->](-3,0)--(12,0);
\draw[thick,->](0,-2)--(0,7);
\draw[<->] (12,0) node[below]{$s$} (0,7) node[left]{$f$}; 
\foreach \t in{-2,-1,...,11}\draw(\t,0.1)--(\t,-0.1);
\foreach \t in{-1,0,...,6}\draw(0.1,\t)--(-0.1,\t);

\draw[fill=black,black](-1,0)circle(0.1);
\draw[fill=black,black](0,1)circle(0.1);
\draw[fill=black,black](1,2)circle(0.1);
\draw[fill=black,black](1,1)circle(0.1);
\draw[fill=black,black](2,2)circle(0.1);
\draw[fill=black,black](3,3)circle(0.1);

\draw[ultra thick,black](-1,0)--(1,2);
\draw[ultra thick,black](1,1)--(1,2);
\draw[ultra thick,black](1,1)--(3,3);

\draw[fill=black,black](7,4)circle(0.1);
\draw[fill=black,black](8,5)circle(0.1);
\draw[fill=black,black](9,6)circle(0.1);
\draw[fill=black,black](9,5)circle(0.1);
\draw[fill=black,black](10,6)circle(0.1);
\draw[fill=black,black](11,7)circle(0.1);

\draw[ultra thick,black](7,4)--(9,6);
\draw[ultra thick,black](9,5)--(9,6);
\draw[ultra thick,black](9,5)--(11,7);

\draw[ultra thick,red,->](1,2)--(2,3);
\draw[ultra thick,red,->](3,3)--(4,4);
\draw[ultra thick,red,->](9,6)--(10,7);
\draw[ultra thick,red,->](11,7)--(12,8);

\end{tikzpicture}

\caption{$[S/h_0,F_0]^{\cA(1)\dual}$}
\end{figure}

This is periodic,where the periodicity shifts degree by $(8,4,4)$. Since $[S/h_0,F_0/h_1^\infty]^{\cA(1)\dual}$ is a colimit, it is essential to know the maps over which we are taking the colimit. First let us take a look at the maps induced by multiplying by $h_1$ (we abbreviate $\Sigma^{-i,-i,-i}$ to $\Sigma^{-i}$ in this diagram):\\

\[
\scalebox{0.87}{
\xymatrix @R=2em {
\ar[r]^-{h_1}&[S/h_0,\Sigma^{-1}F_0]\ar[r] \ar[d]_{h_1\circ\Sigma^{-1}}& [S/h_0,\Sigma^{-1}F_0/h_1]\ar[r]\ar[d]&\Sigma^{2,0,1} [S/h_0,\Sigma^{-1}F_0]\ar[d]^{id}\ar[r]^-{h_1}& \\
\ar[r]^-{h_1^2}&[S/h_0,\Sigma^{-2}F_0]\ar[d]_{h_1\circ\Sigma^{-1}}\ar[r]&[S/h_0,\Sigma^{-2}F_0/h_1^2]\ar[d]\ar[r]&\Sigma^{3,1,2}[S/h_0,\Sigma^{-2}F_0]\ar[d]\ar[r]^-{h_1^2}&\\
 & & & & \\
}
}
\]
\begin{equation}
\label{colim}
\end{equation}

All rows are exact. The colimit of the column on the left is $coker(h_1^k)$, while the colimit of the column on the right is $ker(h_1^k)$. As a result, taking the colimit in the middle would merely be taking the colimits of the cokernel part from the left and the kernel part from the right. There are hidden multiplicative relations between the cokernel and kernel. However, they do not affect the vanishing region, which is our only goal. Here is a more illuminating diagram:\\

\[\scalebox{0.94}{
\xymatrix{
0\ar[r]&coker(h_1)\ar[r] \ar[d]_{h_1\circ\Sigma^{-1}}& [S/h_0,\Sigma^{-1}F_0/h_1]\ar[r]\ar[d]&ker(h_1)\ar[d]^{i}\ar[r]&0\\
0\ar[r]&coker(h_1^2)\ar[d]_{h_1\circ\Sigma^{-1}}\ar[r]&[S/h_0,\Sigma^{-2}F_0/h_1^2]\ar[d]\ar[r]&ker(h_1^2)\ar[d]\ar[r]&0\\
& & & & \\
}}
\]

The maps $i$ on the right column are canonical inclusions, and passing to colimits gives \[\underset{k}{\colim}(coker(h_1^k))\to[S/h_0,F_0/h_1^\infty]\to\underset{k}{\colim}(ker(h_1^k)).\] Working over the dual subalgebra $\cA(1)\dual$ we can calculate $[S/h_0,\Sigma^{-1,1,0}F_0/h_1^\infty]^{\cA(1)\dual}_{*,*,*}$ directly. Furthermore we have:

\begin{prop}\label{inj}
For any $k\in \Z$, $k\geq 1$, the maps $[S/h_0,\Sigma^{-k}F_0/h_1^k]^{\cA(1)\dual}\to [S/h_0,\Sigma^{-k-1}F_0/h_1^{k+1}]^{\cA(1)\dual}$ are injective.
\end{prop}

The result of the calculation is shown in Figure \ref{period8}. The shift in the figure appears as result of Proposition \ref{shift}.

\begin{figure}[h]
\centering
\begin{tikzpicture}[scale=0.6]

\draw[thick,->](-4,0)--(12,0);
\draw[thick,->](0,-2)--(0,7);
\draw[<->] (12,0) node[below]{$s$} (0,7) node[left]{$f$}; 
\foreach \t in{-3,-2,...,11}\draw(\t,0.1)--(\t,-0.1);
\foreach \t in{-1,0,...,6}\draw(0.1,\t)--(-0.1,\t);

\draw[fill=black,black](-1,0)circle(0.1);
\draw[fill=black,black](0,1)circle(0.1);
\draw[fill=black,black](1,2)circle(0.1);
\draw[fill=black,black](1,1)circle(0.1);
\draw[fill=black,black](2,2)circle(0.1);
\draw[fill=black,black](3,3)circle(0.1);

\draw[ultra thick,black](-1,0)--(1,2);
\draw[ultra thick,black](1,1)--(3,3);
\draw[ultra thick,black](1,1)--(1,2);

\draw[fill=black,black](7,4)circle(0.1);
\draw[fill=black,black](8,5)circle(0.1);
\draw[fill=black,black](9,6)circle(0.1);
\draw[fill=black,black](9,5)circle(0.1);
\draw[fill=black,black](10,6)circle(0.1);
\draw[fill=black,black](11,7)circle(0.1);

\draw[ultra thick,black](7,4)--(9,6);
\draw[ultra thick,black](9,5)--(9,6);
\draw[ultra thick,black](9,5)--(11,7);

\draw[fill=black,black](-1,1)circle(0.1);
\draw[fill=black,black](-2,0)circle(0.1);
\draw[fill=black,black](-3,0)circle(0.1);
\draw[fill=black,black](-4,-1)circle(0.1);
\draw[ultra thick,pink,->](-1,1)--(-3,-1);
\draw[ultra thick,pink,->](-3,0)--(-5,-2);
\draw[fill=black,black](7,5)circle(0.1);
\draw[fill=black,black](6,4)circle(0.1);
\draw[fill=black,black](5,4)circle(0.1);
\draw[fill=black,black](4,3)circle(0.1);
\draw[ultra thick,pink,->](7,5)--(5,3);
\draw[ultra thick,pink,->](5,4)--(3,2);

\end{tikzpicture}

\caption{$[S/h_0,\Sigma^{-1,1,0}F_0/h_1^\infty]^{\cA(1)\dual}_{*,*,*}$}\label{period8}
\end{figure}

This is periodic, with a periodicity degree shift of $(8,4,4)$, just as with $[S/h_0,F_0]^{\cA(1)\dual}$. Note that $[S/h_0,\Sigma^{-1,1,0}F_0/h_1^\infty]^{\cA(1)\dual}_{*,*,*}$ differs from the classical $[S/h_0,S]^{\cA_{cl}(1)\dual}_{*,*}$ with two extra negative $h_1$-towers associated to each "lighting flash". The element in degree $(-1,0,-1)$ in the first pattern is generated by $\tau$ with a shift.

Recall the self-map $P$ on $S/h_0$ acts injectively as can be seen in Figure \ref{period8}. Combining this with the long exact sequence:
\[
\xymatrix @R=1.5ex{\cdots\ar[r]& [S/(h_0,P),F_0/h_1^\infty]_{s,f,w}^{\cA(1)\dual}\ar[r]& [S/h_0,F_0/h_1^\infty]_{s,f,w}^{\cA(1)\dual}\ar[r]^-{P}& \\
 \ar[r]^-{P} &[S/h_0, F_0/h_1^\infty]_{s+8,f+4,w+4}^{\cA(1)\dual}\ar[r] &[S/(h_0,P),F_0/h_1^\infty]_{s-1,f+1,w}^{\cA(1)\dual}\ar[r]& \cdots\\
}
\]
\noindent gives $[S/(h_0,P),\Sigma^{-1,1,0}F_0/h_1^\infty]_{*,*,*}^{\cA(1)\dual}$ as in Figure \ref{figure2}.

\begin{rmk}\label{also inj}
Analogously to Proposition \ref{inj}, for any $k\in \Z$, $k\geq 1$, the following maps are also injective: \[[S/(h_0,P),\Sigma^{-k}F_0/h_1^k]^{\cA(1)\dual}\to [S/(h_0,P),\Sigma^{-k-1}F_0/h_1^{k+1}]^{\cA(1)\dual}.\] 
\end{rmk}

\begin{figure}[h]
\centering
\begin{tikzpicture}[scale=0.6]

\draw[thick,->](-12,0)--(2,0);
\draw[thick,->](0,-5)--(0,4);
\draw[<->] (2,0) node[below]{$s$} (0,4) node[left]{$f$}; 
\foreach \t in{-11,-10,...,1}\draw(\t,0.1)--(\t,-0.1);
\foreach \t in{-4,-3,...,3}\draw(0.1,\t)--(-0.1,\t);

\draw[fill=black,black](-10,-3)circle(0.1);
\draw[fill=black,black](-9,-2)circle(0.1);
\draw[fill=black,black](-8,-1)circle(0.1);
\draw[fill=black,black](-8,-2)circle(0.1);
\draw[fill=black,black](-7,-1)circle(0.1);
\draw[fill=black,black](-6,0)circle(0.1);

\draw[ultra thick,black](-10,-3)--(-8,-1);
\draw[ultra thick,black](-8,-2)--(-8,-1);
\draw[ultra thick,black](-8,-2)--(-6,0);

\draw[fill=black,black](-10,-2)circle(0.1);
\draw[fill=black,black](-11,-3)circle(0.1);
\draw[fill=black,black](-12,-3)circle(0.1);
\draw[fill=black,black](-13,-4)circle(0.1);
\draw[ultra thick,pink,->](-10,-2)--(-12,-4);
\draw[ultra thick,pink,->](-12,-3)--(-14,-5);

\end{tikzpicture}

\caption{$[S/(h_0,P),\Sigma^{-1,1,0}F_0/h_1^\infty]^{\cA(1)\dual}_{*,*,*}$}\label{figure2}
\end{figure}

Next we will use the Cartan-Eilenberg spectral sequence to bootstrap our result from $\cA(1)\dual$-homology to $\cA\dual$-homology. The Cartan-Eilenberg spectral sequence converges when the input is a bounded-below $\cA\dual$-comodule. We will obtain a vanishing region for each finite stage $[S/(h_0,P),\Sigma^{-k}F_0/h_1^k]^{\cA\dual}$ and then deduce the vanishing region for $[S/(h_0,P),F_0/h_1^\infty]^{\cA\dual}$ by passing to the colimit.
\[
\scalebox{0.92}{
\xymatrix{
[S/(h_0,P),\Sigma^{-1}F_0/h_1]^{\cA(1)\dual}\ar[r]\ar@{~>}[d]&[S/(h_0,P),\Sigma^{-2}F_0/h_1^2]^{\cA(1)\dual}\ar[r]\ar@{~>}[d]_{CESS}&\cdots\ar[r]&[S/(h_0,P),F_0/h_1^\infty]^{\cA(1)\dual}\\
[S/(h_0,P),\Sigma^{-1}F_0/h_1]^{\cA\dual}\ar[r]&[S/(h_0,P),\Sigma^{-2}F_0/h_1^2]^{\cA\dual}\ar[r]&\cdots\ar[r]&[S/(h_0,P),F_0/h_1^\infty]^{\cA\dual}
}
}\\
\]

Going from $\cA(1)\dual$ to $\cA\dual$ is too big of a step, so we first calculate $[S/(h_0,P),F_0/h_1^\infty]^{\cA(2)\dual}$, where
\[\cA(2)\dual=\M_2[\tau_0,\tau_1,\tau_2,\xi_1,\xi_2]/(\tau_0^2=\tau\xi_1,\tau_1^2=\tau\xi_2,\tau_2^2,\xi_1^4,\xi_2^2).\] 
To do this, we will use a sequence of normal maps of Hopf algebras:
\[\cA(2)\dual\to \cA(2)\dual/\xi_1^2\to \cA(2)\dual/(\xi_1^2,\xi_2)\to \cA(2)\dual/(\xi_1^2,\xi_2,\tau_2)=\cA(1)\dual.\]

First we consider the Cartan-Eilenberg spectral sequence corresponding to the extension $$E(\tau_2)\to \cA(2)\dual/(\xi_1^2,\xi_2)\to \cA(1)\dual.$$ The element $\tau_2$, which has degree $(6,1,3)$, corresponds to $h_{30}$ in the May spectral sequence. The $\cA(1)\dual$-coaction on $E(\tau_2)$ is trivial for degree reasons. So we start with the $E_1=E_2$-page, and deduce a vanishing region on $[S/(h_0,P),F_0/h_1^\infty]^{\cA(2)\dual/(\xi_1^2,\xi_2)}$.
\[
\xymatrix{
[S/(h_0,P),\Sigma^{-1}F_0/h_1]^{\cA(1)\dual}\otimes\M_2[h_{30}]\ar[r]\ar@{=>}[d]&\cdots\ar[r]&[S/(h_0,P),F_0/h_1^\infty]^{\cA(1)\dual}\otimes\M_2[h_{30}]\\
[S/(h_0,P),\Sigma^{-1}F_0/h_1]^{\cA(2)\dual/(\xi_1^2,\xi_2)}\ar[r]&\cdots\ar[r]&[S/(h_0,P),F_0/h_1^\infty]^{\cA(2)\dual/(\xi_1^2,\xi_2)}
}
\]

For the normal extension $E(\beta)\to \Gamma \to C$ of Hopf algebras we state a motivic version of \cite[Lemma 4.10]{Kra}, which gives a relationship between the vanishing region for $[M,N]^\Gamma$ and the vanishing condition of $[M,N]^C$ together with the two "slopes" associated to $\beta$. Note that if $\beta$ has degree $(s_0,f_0,w_0)$, then $\frac{f_0}{s_0}$ and $\frac{f_0}{w_0}$ are the slopes of the projections of $(s_0,f_0,w_0)$ onto the plane $w=0$ and the plane $s=0$.

\begin{thm}\label{KraLem410}
Let $E(\alpha)\to \Gamma \xrightarrow{q} C$ be a normal extension of Hopf algebras and $M$,$N\in Stab(\Gamma)$. Suppose $\beta$ is an element in $[S,S]^E$ of degree $(s_0,f_0,w_0)$ with $s_0,f_0,w_0$ all positive. Its image in $[S,S]^\Gamma$ (which we also call $\beta$) acts on $[M,N]^\Gamma$. Suppose for some $a,b,c,m,c_0\in \R$ with $a,b>0$ and $m\geq\frac{f_0}{s_0}>0$, the group $[q_*(M),q_*(N)]^C$ vanishes when $f>as+bw+c$ and also vanishes when $f>ms+c_0$. Then
\begin{enumerate}
    \item if $f_0\leq as_0+bw_0$, or $\beta$ acts nilpotently on $[M,N]^\Gamma$, then $[M,N]^\Gamma$ has a parallel vanishing region. In other words, it vanishes when $f>as+bw+c'$ for some constant $c'$ and also vanishes when $f>ms+c_0$.
    \item otherwise, $[M,N]^\Gamma$ vanishes when $f>\frac{mbw_0-f_0(m-a)}{bw_0-s_0(m-a)}s+\frac{bf_0-mbs_0}{bw_0-s_0(m-a)}w+c'$ and vanishes when $f>ms+c_0$.
\end{enumerate}
\end{thm}

\begin{rmk}
The additional vanishing plane $f>ms+c_0$ generalizes the bounded below condition. In the classical setting, we have that $[M,N]^\Gamma$ vanishes when $s<c_0$, but due to the negative $h_1$-towers we do not have a vertical vanishing plane. So we adjust the "$\infty$-slope" plane to be $f=ms+c_0$ to fulfill our purpose. This bound does not affect the periodicity region we study here, so we omit it henceforth.
\end{rmk}

\begin{pf}[Proof of Theorem~\ref{KraLem410}]
If $\beta$ has $f_0\leq as_0+bw_0$, then $\beta$ multiples of classes in $[M,N]^C$ will lie under the existing vanishing planes.

If $f_0> as_0+bw_0$, then every infinite $\beta$ tower will contain classes lying outside of the rigion $f>as+bw+c$. If $\beta$ acts nilpotently, however, then there is a micimum length for all $\beta$-towers, and so we can still get a parallel vanishing plane $f>as+bw+c'$ on $[M,N]^\Gamma$ by adjusting the $f$-intercept.

Now we turn to case $(2)$. If $f_0> as_0+bw_0$ and $\beta$ acts non-nilpotently, then there must exist an element $x\in [M,N]^\Gamma$ for which the classes $\beta^kx$ are not zero on the $E_\infty$ page of the Cartan-Eilenberg spectral sequence for every $k$. Thus no matter how we move up the existing vanishing plane $f>as+bw+c$, some $\beta$ multiples of $x$ will lie above the plane. Instead, we will find a new vanishing plane $f>a's+b'w+c'$ for $a',b',c'\in \R$. The new vanishing region $f>a's+b'w+c'$ must satisfy the condition $f_0\leq a's_0+b'w_0+c'$. This plane is spanned by the direction of $\beta$ and the intersecting line of the two existing vanishing planes. Hence we can solve to obtain $a'=\frac{mbw_0-f_0(m-a)}{bw_0-s_0(m-a)}$ and $b'=\frac{bf_0-mbs_0}{bw_0-s_0(m-a)}$.
\end{pf}

\begin{rmk}\label{handy}
In the relevant cases, the starting vanishing regions will have $b=0$. One can think of these as 2-dimensional cases stated in 3-dimensional language.\\ We rewrite the conditions and the results of Theorem \ref{KraLem410} as the following:
Suppose for some $a,c,m,c_0\in \R$ with $a>0$ and $m\geq\frac{f_0}{s_0}>0$, the group $[q_*(M),q_*(N)]^C$ vanishes when $f>as+c$ and also vanishes when $f>ms+c_0$. Then:
\begin{enumerate}
    \item if $f_0\leq as_0$, or $\beta$ acts nilpotently on $[M,N]^\Gamma$, then $[M,N]^\Gamma$ has a parallel vanishing region. That is to say, it vanishes when $f>as+c'$ for some constant $c'$, and also vanishes when $f>ms+c_0$,
    \item if otherwise, then $[M,N]^\Gamma$ vanishes when $f>\frac{f_0}{s_0}s+c'$ for some constant $c'$, and vanishes when $f>ms+c_0$.
\end{enumerate}
\end{rmk}

\begin{rmk}
Similarly, we could generalize to the statement that $[q_*(M),q_*(N)]^C$ vanishes when $f>\varphi(s,w)$ where $\varphi:\R^2\to\R$ is a smooth function. Then the gradient $v(-,-)=(\frac{\bdry\varphi}{\bdry s}(-),\frac{\bdry\varphi}{\bdry w}(-))$ is a vector field. Now we would like to consider $g=\underset{(s_0,w_0)}{Min}|v(s_0,w_0)|$ and compare $g$ with $\frac{f_0}{s_0}$ and $\frac{f_0}{w_0}$. The conditions can be rewritten as follows:
\begin{enumerate}
    \item if $\frac{f_0}{s_0}\leq g$ or $\frac{f_0}{w_0}\leq g$, or $\beta$ acts nilpotently, then $[M,N]^\Gamma$ has a parallel vanishing region.
    \item if both $\frac{f_0}{s_0}$ and $\frac{f_0}{w_0}>g$, and $\beta$ acts non-nilpotently, then we must modify the vanishing region of $[M,N]^\Gamma$. However, it takes some work to write down a precise modification, so we omit it here.
\end{enumerate}
\end{rmk}

\begin{rmk}
From the cofiber sequence $S\xrightarrow{h_0^k} S\to S/h_0^k$ we can take tensor duals to derive the fiber sequence $D(S/h_0^k)\to S\to S$. Since $D(S/h_0^k) \simeq \Sigma^{-1,1-k,0}S/h_0^k$, we have
\[[S/h_0^k, S]_{s,f,w}=[S, D(S/h_0^k)]_{s,f,w}=[S, S/h_0^k]_{s+1,f+k-1,w.}\]
Because $S/h_0^k$ is compact in $Stab(\cA\dual)$, smashing the second slot with some $N\in Stab(\cA\dual)$, we get \[[S/h_0^k,N]_{s,f,w}\iso[S,D(S/h_0^k)\smsh N]_{s,f,w}\iso [S,S/h_0^k\smsh N]_{s+1,f+k-1,w.}\]
As a result $\beta\in [S,S]^\Gamma$ acts on $[M,N]^\Gamma$ for compact $M\in Stab(\cA\dual)$, since $\beta$ acts on $[S,DM\smsh N]^\Gamma$.
\end{rmk}

The group $[S/(h_0,P),\Sigma^{-1,1,0}F_0/h_1^\infty]_{*,*,*}^{\cA(1)\dual}$ has a single "lighting flash" pattern along with two negative $h_1$-towers (see Figure 
\ref{figure2}), so the vanishing region to start off with is $f>c$ (We obtain the same vanishing region of $[S/(h_0,P),\Sigma^{-1,1,0}(\Sigma^{-k}F_0/h_1^k)]_{*,*,*}^{\cA(1)\dual}$ for each $k$, since the maps we are taking colimit over are injections by Remark \ref{also inj}.) In our case, $[M,N]^\Gamma=[S/(h_0,P),\Sigma^{-1,1,0}F_0/h_1^\infty]_{*,*,*}^{\cA(1)\dual}$, and we will apply Theorem \ref{KraLem410} in the following three cases: (i) $\beta$ is $\tau_2$ of degree $(6,1,3)$; (ii) $\beta$ is $\xi_2$ of degree $(5,1,3)$; (iii) $\beta$ is $\xi_1^2$ of degree $(3,1,2)$.

Recall that we are working with the Cartan-Eilenberg spectral sequence \[[S/(h_0,P),\Sigma^{-1,1,0}(\Sigma^{-k}F_0/h_1^k)]^{\cA(1)\dual}\otimes\M_2[h_{30}]\Rightarrow [S/(h_0,P),\Sigma^{-1,1,0}(\Sigma^{-k}F_0/h_1^k)]^{\cA(2)\dual/(\xi_1^2,\xi_2)}.\]
There cannot be any differentials for degree reasons. By Theorem \ref{KraLem410} the element $h_{30}$ will bring us a vanishing region $f>\frac{1}{6}s+c_1$ for each $k$, where $c_1$ is some constant (we obtain the same constant for all $k$). Passing to the colimit, we conclude that $[S/(h_0,P),\Sigma^{-1,1,0}F_0/h_1^
\infty]^{\cA(2)\dual/(\xi_1^2,\xi_2)}$ shares the same vanishing region $f>\frac{1}{6}s+c_1$.\\

The second step is to consider the normal extension in which we add $\xi_2$, corresponding to the class $h_{21}$: 
\[E(\xi_2)\to \cA(2)\dual/\xi_1^2\to \cA(2)\dual/(\xi_1^2,\xi_2).\]
The $\cA(2)\dual/(\xi_1^2,\xi_2)$-coaction on $E(\xi_2)$ is trivial. We have $E_2$-pages as the first row:
\[\scalebox{0.92}{
\xymatrix{
[S/(h_0,P),\Sigma^{-1}F_0/h_1]^{\cA(2)\dual/(\xi_1^2,\xi_2)}\otimes\M_2[h_{21}]\ar[r]\ar@{=>}[d]&\cdots\ar[r]&[S/(h_0,P),F_0/h_1^\infty]^{\cA(2)\dual/(\xi_1^2,\xi_2)}\otimes\M_2[h_{21}]\\
[S/(h_0,P),\Sigma^{-1}F_0/h_1]^{\cA(2)\dual/\xi_1^2}\ar[r]&\cdots\ar[r]&[S/(h_0,P),F_0/h_1^\infty]^{\cA(2)\dual/\xi_1^2}
}}
\]

The spectral sequence collapses at the $E_2$-page. This is because in the May spectral sequence over $\cA(2)$ or $\cA$, there is a differential $d_1(h_{30})=h_1h_{21}+h_2h_{20}$, but $h_0$ and $h_2$ are zero in the group $[S/(h_0,P),\Sigma^{-1,1,0}(\Sigma^{-k}F_0/h_1^k)]^{\cA(2)\dual/(\xi_1^2,\xi_2)}$. As a result, $h_{21}$ is also non-nilpotent. For some constant $c_2$, the vanishing region of $[S/(h_0,P),\Sigma^{-1,1,0}(\Sigma^{-k}F_0/h_1^k)]^{\cA(2)\dual/\xi_1^2}$ is $f>\frac{1}{5}s+c_2$ for each $k$ according to Theorem \ref{KraLem410}, and the same is true for the colimit $[S/(h_0,P),\Sigma^{-1,1,0}F_0/h_1^\infty]^{\cA(2)\dual/\xi_1^2}$.\\

Next we consider the Cartan-Eilenberg spectral sequence corresponding to the extension: \[E(\xi_1^2)\to \cA(2)\dual\to \cA(2)\dual/\xi_{1.}^2\]
Here the class $\xi_1^2$ corresponds to the class $h_2$ in the May spectral sequence. The $\cA(2)\dual/\xi_1^2$-coaction on $E(\xi_1^2)$ is trivial as well. We have $E_2$-pages as in the first row:

\[\xymatrix{
[S/(h_0,P),\Sigma\inv F_0/h_1]^{\cA(2)\dual/\xi_1^2}\otimes\M_2[h_2]\ar[r]\ar@{=>}[d]&\cdots\ar[r]&[S/(h_0,P),F_0/h_1^\infty]^{\cA(2)\dual/\xi_1^2}\otimes\M_2[h_2]\\
[S/(h_0,P),\Sigma\inv F_0/h_1]^{\cA(2)\dual}\ar[r]&\cdots\ar[r]&[S/(h_0,P),F_0/h_1^\infty]^{\cA(2)\dual}
}
\]

We do have some non-zero differentials appear. In the previous steps, by introducing $[\tau_2]=(6,1,3)$ and $[\xi_2]=(5,1,3)$, which give rise to non-nilpotent elements in $\Ext$, we arrived a vanishing region of $f>\frac{1}{5}s+c_3$, where $c_3$ is a constant. However $[\xi_1^2]=(3,1,2)$ is nilpotent since $h_2^4=0$ in $\Ext_{\cA(2)\dual}$ and $\Ext$. 

Moving from $\cA(2)\dual$ to $\cA\dual$, we have many more elements to introduce. However those elements won't satisfy $\frac{f}{s}>\frac{1}{5}$. By Theorem \ref{KraLem410} (or Remark~\ref{handy}), for each $k$, $[S/(h_0,P),\Sigma^{-1,1,0}(\Sigma^{-k}F_0/h_1^k)]^\cA$ vanishes if $f=\frac{1}{5}s+c_3$. Since the vanishing plane passes through the point $(-6,0,-1)+3\cdot(3,1,2)=(3,3,5)$, the constant $c_3$ is $\frac{12}{5}$ and the region $f>\frac{1}{5}s+\frac{12}{5}$ would be carried through to $\cA\dual$. We conclude that
\begin{prop}\label{vrk=1}
The group $[S/(h_0,P),\Sigma^{-1,1,0}F_0/h_1^\infty]_{s,f,w}$ has a vanishing region of $f>\frac{1}{5}s+\frac{12}{5}$.
\end{prop}

Note that it is possible for many reasons that the vanishing region we have found is not optimal. First, we could consider the "slope" of the motivic weight side $\frac{f}{w}$ instead of $\frac{f}{s}$ under certain bounded below conditions. Second, if other elements were included, more differentials would occur, allowing for a larger vanishing region. More calculation is required to clarify these cases.

\section{The motivic periodicity theorem}

Let $F_0$ and $F_{01}$ still be the same as in Definition \ref{F0}, so that $[S,\Sigma^{-1,1,0}F_0/h_1^\infty]_{s,f,w}\iso [S,F_{01}]_{s,f,w}$ when $f>\frac{1}{2}s+3$. Given a self-map $\theta$ on $S/h_0^k$ let us recall the diagram where the first row is exact:\\
\[
\scalebox{0.78}{
\xymatrix{
[S/(h_0^k,\theta),\Sigma^{-1,1,0}F_0/h_1^\infty]\ar[r]& [S/h_0^k,\Sigma^{-1,1,0}F_0/h_1^\infty]\ar[r]^-{\theta} \ar[d] & [S/h_0^k,\Sigma^{-1,1,0}F_0/h_1^\infty]\ar[d]\ar[r]& \Sigma^{-1,1,0}[S/(h_0^k,\theta),\Sigma^{-1,1,0}F_0/h_1^\infty] \\
 & [S,\Sigma^{-1,1,0}F_0/h_1^\infty]\ar[r]^-{P_r(-)} & [S,\Sigma^{-1,1,0}F_0/h_1^\infty] &
}
}
\]\\

The vertical maps are isomorphisms whenever $f>\frac{1}{2}s+\frac{3}{2}-k$ due to Corollary \ref{ud}. We would like to further restrict the condition to $f>\frac{1}{2}s+3-k$ in order to eliminate the indeterminacy. The vanishing condition on $[S/(h_0^k,\theta),\Sigma^{-1,1,0}F_0/h_1^\infty]$, which is the same as the vanishing condition on $[S/(h_0^k,\theta),F_{01}]_{s,f,w}$, tells us whether $\theta$ is an isomorphism.\\

In the previous section, we established the case when $k=1$, given in Proposition \ref{vrk=1}. We show in Figure \ref{periodmore} the $(2^{r+1},2^r,2^r)$-periodic pattern for $[S/h_0^k,\Sigma^{-1,1,0}F_0/h_1^\infty]^{\cA(1)\dual}$, where $k\leq 2^r$. By an analogous computation,  one can see that for a general positive integer $k\leq 2^r$, the groups $[S/(h_0^k,P^{2^{r-2}}),F_{01}]_{s,f,w}$ admit a parallel vanishing region as in the $k=1$ case.

\begin{figure}[!h]
\centering
\begin{tikzpicture}[scale=0.6]

\draw[thick,->](-4,0)--(12,0);
\draw[thick,->](0,-2)--(0,10);
\draw[<->] (12,0) node[below]{$s$} (0,10) node[left]{$f$};

\draw[fill=black,black](-1,0)circle(0.1);
\draw(-1,0) node[below]{$(-1,0,-1)$};
\draw[fill=black,black](0,1)circle(0.1);
\draw[fill=black,black](1,2)circle(0.1);
\draw[fill=black,black](1,4)circle(0.1);
\draw[fill=black,black](2,5)circle(0.1);
\draw[fill=black,black](3,6)circle(0.1);
\draw(3,6) node[above]{$(3,k+1,1)$};

\draw(4,3) node[below]{$\cdots$};
\draw(1,3) node[left]{$\vdots$};
\draw(11,7) node[left]{$\vdots$};

\draw[ultra thick,black](-1,0)--(1,2);
\draw[ultra thick,black](1,4)--(3,6);

\draw[fill=black,black](9,4)circle(0.1);
\draw[fill=black,black](10,5)circle(0.1);
\draw[fill=black,black](11,6)circle(0.1);
\draw[fill=black,black](11,8)circle(0.1);
\draw[fill=black,black](12,9)circle(0.1);
\draw[fill=black,black](13,10)circle(0.1);
\draw(13,10) node[above]{$(2^{r+1}+3,2^r+k+1,2^r+2-1)$};

\draw[ultra thick,black](9,4)--(11,6);
\draw[ultra thick,black](11,8)--(13,10);

\draw[fill=black,black](-1,4)circle(0.1);
\draw[fill=black,black](-2,3)circle(0.1);
\draw[fill=black,black](-3,0)circle(0.1);
\draw[fill=black,black](-4,-1)circle(0.1);
\draw[ultra thick,pink,->](-1,4)--(-3,2);
\draw[ultra thick,pink,->](-3,0)--(-5,-2);
\draw[fill=black,black](9,8)circle(0.1);
\draw[fill=black,black](8,7)circle(0.1);
\draw[fill=black,black](7,4)circle(0.1);

\draw[fill=black,black](6,3)circle(0.1);
\draw[ultra thick,pink,->](9,8)--(7,6);
\draw[ultra thick,pink,->](7,4)--(5,2);

\end{tikzpicture}
\caption{$[S/h_0^k,\Sigma^{-1,1,0}F_0/h_1^\infty]^{\cA(1)\dual}_{*,*,*}$}\label{periodmore}
\end{figure}

We have the following lemma for the $f$-intercept:

\begin{lemma}[Corollary of {\cite[Lemma 5.4]{Kra}}]\label{Kra5.4}
Let $M,N\in Stable(\cA\dual)$ with $M$ compact. Let $M_1=M/\theta_1$ be the cofiber of the self-maps $\Sigma^{s_1,f_1,w_1}M\xrightarrow{\theta_1}M$, and let $M_2=M/(\theta_1,\theta_2)$ be the cofiber of the self-maps $\Sigma^{s_2,f_2,w_2}M/\theta_1\xrightarrow{\theta_2}M/\theta_1$. Define $M_1'$ and $M_2'$ with respect to the self-maps $\Sigma^{s_1',f_1',w_1'}M\xrightarrow{\theta_1'}M$ and $\Sigma^{s_2',f_2',w_2'}M/\theta_1'\xrightarrow{\theta_2'}M/\theta_1'$ in the same way. Suppose $\theta_i$ and $\theta_i'$ are parallel, i.e. $(s_i,f_i,w_i)=\lambda_i(s_i',f_i',w_i')$ where $\lambda_i$ are non-zero real numbers and $i=1,2$. 

Further let $a,b\in \R$ and suppose $f_i>as_i+bw_i$ and $f_i'>as_i'+bw_i'$ for $i=1,2$. We make the convention that the $f$-intercept is $\infty$ if there is no such vanishing plane. Then the minimal $f$-intercepts of the vanishing planes parallel to $f=as+bw$ on $[M_2,N]$ and $[M_2',N]$ agree.
\end{lemma}

\begin{pf}[Proof of Lemma \ref{Kra5.4}]
We construct the iterated cofiber $L_1=M/(\theta_1,\theta_1')$ and $L_2=M/(\theta_1,\theta_2,\theta_1',\theta_2')$. Since $f_i>as_i+bw_i$ and $f_i'>as_i'+bw_i'$ for $i=1,2$, the minimal $f$-intercepts for the vanishing planes parellel to $f=as+bw$ agree on $[M_i,N]$, $[M_i',N]$ and $[L_i,N]$ by inductively applying Lemma \ref{lemma 5.2}.

Note that the notation for $L_1$ and $L_2$ is ambiguous. The notation does not indicate that $M/\theta_1$ should admit a $\theta_1'$ self-map or vice versa. Because of the uniqueness of (homological) self-maps that Krause has shown in \cite[Sec. 4]{Kra}, there is a self-map $\theta_1''$ compatible with both $\theta_1$ and $\theta_1'$, which acts on $M$ by a power of $\theta_1$, and by a power of $\theta_2$. We will take $L_1$ to be the cofiber of the self-map $\theta_1''$. Similarly, there exists a self-map $\theta_2''$ on $L_1$ that acts on $M_1$ by a power of $\theta_2$, and on $M_1'$ by a power of $\theta_2'$. So we can set $L_2$ as the cofiber of the self-map $\theta_2''$.
\end{pf}

\begin{rmk}
Krause's proof of the uniqueness of self-maps is in the classical setting, yet for the $\C$-motivic case the proof is analogous.
\end{rmk}

\begin{rmk}
The cofiber sequences arising from the Verdier's axiom and the $3\times 3$ lemma offer an alternative way to issue the vanishing condition of $[S/(h_0^k,P^{2^{r-2}}),F_{01}]_{s,f,w}$. Let $m,n,l,l'\in \N$ be positive with $m\leq 4l$ and $m+n\leq 4(l+l')$. We have the following cofiber sequences:
\[S/h_0^m\to S/h_0^{m+n}\to S/h_0^n\]
\[S/(h_0^m,P^{l+l'})\to S/(h_0^{m+n},P^{l+l'})\to S/(h_0^n,P^{l+l'})\]
\[S/(h_0^m,P^l)\to S/(h_0^m,P^{l+l'})\to S/(h_0^m,P^{l'}).\]
Passing to the induced long exact sequences in homology, we conclude that for $k\leq 2^r$, the groups $[S/(h_0^k,P^{2^{r-2}}),F_{01}]_{s,f,w}$ admit the same vanishing condition as $[S/(h_0,P),F_{01}]_{s,f,w}$.
\end{rmk}

It follows that for any $k\leq 2^r$ and any self-map $\theta=P^{2^{r-2}}$ of $S/h_0^k$, the corresponding groups $[S/(h_0^k,\theta), F_{01}]$ have a vanishing region of $f>\frac{1}{5}s+\frac{12}{5}$. Combining with Theorem~\ref{selfmap}, we arrive at the motivic version of Theorem~\ref{cpt}:

\begin{thm}[Another way of stating Theorem~\ref{mpt}]
For $r\geq 2$, the Massey product operation $P_r(-):=\langle h_{r+1}, h_0^{2^r},-\rangle$ is uniquely defined on $\Ext^{s,f,w}=H^{s,f,w}(\cA)$ when $s>0$ and $f>\frac{1}{2}s+3-2^r$.

Furthermore, for $f>\frac{1}{5}s+\frac{12}{5}$, 
\[P_r\colon [S,F_{01}]_{s,f,w}\xrightarrow{P_r(-)}[S,F_{01}]_{s+2^{r+1},f+2^r,w+2^r.}\]
is an isomorphism between $h_0$ and $h_1$-torsions.
\end{thm}

\end{document}